\documentclass{article}
\usepackage{amssymb,amsmath,euscript}

\begin{document}

\def\Sp{\mathrm {Sp}}
\def\U{\mathrm U}
\def\SOS{\mathrm {SO}^*}

\def\R{{\mathbb R}}
\def\C{{\mathbb C}}
\def\T{{\mathbb T}}
\def\Z{{\mathbb Z}}

\def\B{{\rm B}}
\def\ov{\overline}
\def\wt{\widetilde}

\def\phi{\varphi}
\def\kappa{\varkappa}
\def\epsilon{\varepsilon}

\newcommand{\Vir}{\mathop{\mathrm{Vir}}\nolimits}

\newcommand{\Vect}{\mathop{\mathrm{Vect}}\nolimits}


\renewcommand{\theequation}{\arabic{section}.\arabic{equation}}
\newcounter{punct}[section]


\renewcommand{\thepunct}{\thesection.\arabic{punct}}
\def\punct{\refstepcounter{punct}{\arabic{section}.\arabic{punct}.  }}


\def\SS{\smallskip}

\newtheorem{theorem}{Theorem}[section]
\newtheorem{proposition}[theorem]{Proposition}
\newtheorem{prop}[theorem]{Proposition}
\newtheorem{lemma}[theorem]{Lemma}
\newtheorem{cor}[theorem]{Corollary}
\newtheorem{corollary}[theorem]{Corollary}
\newtheorem{observation}[theorem]{Observation}


\def\cA{\mathcal A}
\def\cB{\mathcal B}
\def\cC{\mathcal C}
\def\cD{\mathcal D}
\def\cE{\mathcal E}
\def\cF{\mathcal F}
\def\cG{\mathcal G}
\def\cH{\mathcal H}
\def\cJ{\mathcal J}
\def\cI{\mathcal I}
\def\cK{\mathcal K}
\def\cL{\mathcal L}
\def\cM{\mathcal M}
\def\cN{\mathcal N}
\def\cO{\mathcal O}
\def\cP{\mathcal P}
\def\cQ{\mathcal Q}
\def\cR{\mathcal R}
\def\cS{\mathcal S}
\def\cT{\mathcal T}
\def\cU{\mathcal U}
\def\cV{\mathcal V}
\def\cW{\mathcal W}
\def\cX{\mathcal X}
\def\cY{\mathcal Y}
\def\cZ{\mathcal Z}

 \def\cBS{\mathcal {BS}} 


\def\bfm{\mathbf m}


\def\frA{\mathfrak A}
\def\frB{\mathfrak B}
\def\frC{\mathfrak C}
\def\frD{\mathfrak D}
\def\frE{\mathfrak E}
\def\frF{\mathfrak F}
\def\frG{\mathfrak G}
\def\frH{\mathfrak H}
\def\frJ{\mathfrak J}
\def\frK{\mathfrak K}
\def\frL{\mathfrak L}
\def\frM{\mathfrak M}
\def\frN{\mathfrak N}
\def\frO{\mathfrak O}
\def\frP{\mathfrak P}
\def\frQ{\mathfrak Q}
\def\frR{\mathfrak R}
\def\frS{\mathfrak S}
\def\frT{\mathfrak T}
\def\frU{\mathfrak U}
\def\frV{\mathfrak V}
\def\frW{\mathfrak W}
\def\frX{\mathfrak X}
\def\frY{\mathfrak Y}
\def\frZ{\mathfrak Z}

\def\fra{\mathfrak a}
\def\frb{\mathfrak b}
\def\frc{\mathfrak c}
\def\frd{\mathfrak d}
\def\fre{\mathfrak e}
\def\frf{\mathfrak f}
\def\frg{\mathfrak g}
\def\frh{\mathfrak h}
\def\fri{\mathfrak i}
\def\frj{\mathfrak j}
\def\frk{\mathfrak k}
\def\frl{\mathfrak l}
\def\frm{\mathfrak m}
\def\frn{\mathfrak n}
\def\fro{\mathfrak o}
\def\frp{\mathfrak p}
\def\frq{\mathfrak q}
\def\frr{\mathfrak r}
\def\frs{\mathfrak s}
\def\frt{\mathfrak t}
\def\fru{\mathfrak u}
\def\frv{\mathfrak v}
\def\frw{\mathfrak w}
\def\frx{\mathfrak x}
\def\fry{\mathfrak y}
\def\frz{\mathfrak z}

\def\arr{\rightrightarrows}
\def\lr{\leftrightarrow}
\def\lra{\longrightarrow}
\def\llr{\longleftrightarrow}

\def\LR{\Longleftrightarrow}
\def\RA{\Longrightarrow}
\def\LA{\Longleftarrow}

\def\Lr{\Leftrightarrow}
\def\Ra{\Rightarrow}
\def\La{\Leftarrow}

\def\QED{\hfill$\square$}

\def\wh{\widehat}

\def\P{\cP}
\def\Q{\cQ}

\def\SS{\smallskip}

\begin{center}
\bf\Large
On action of the Virasoro algebra on the space of univalent functions

\bigskip

\sc\large
Helene Airault, Yuri A. Neretin%
\footnote{Supported by
the grant FWF, project P19064, Russian Federal Agency for Nuclear Energy,
  NWO.047.017.015 and grant JSPS-RFBR-07.01.91209 }
\end{center}

\bigskip

{\small We obtain explicit expressions for differential
operators defining  the action
of the Virasoro algebra on the space of 
univalent functions. We also obtain an explicit
Taylor decomposition for Schwarzian derivative and 
a formula for the Grunsky coefficients.}

\bigskip


{\bf\punct Virasoro algebra.}
Consider the Lie algebra $\Vect(S^1)$ of vector fields on the circle
$|z|=1$, its basis $L_n:=z^{n+1}\frac d{dz}$
 is numerated by integers, and
$$
[L_n,L_m]=(m-n)L_{n+m}
$$

Recall that the {\it Virasoro algebra} $\Vir$ is the
extended Lie algebra $\Vect(S^1)$; its generators are 
 $L_n$, where $n$ ranges in $\Z$,
and $\zeta$; the commutation relations are
\begin{align}
[L_n,L_m]&=(m-n)L_{m+n}+\frac 1{12} (n^3-n) \delta_{m+n,0}\zeta
\label{vir:1}
\\
[L_n,\zeta]&=0
\nonumber
\end{align} 


{\bf\punct Space of univalent functions $\frK$.}
Denote by   $\frK$  the space of all the functions
$$
f(z)=z+c_1z^2+c_2z^3+\dots
$$
that are univalent in the disk $|z|<1$.
Recall that a function $f$ is univalent,
if $z\ne u$ implies $f(z)\ne f(u)$.
The standard references are \cite{Gol}, \cite{Dur}.

\smallskip


{\bf \punct Action of Virasoro algebra on $\frK$.}
According to Kirillov%
\footnote{In some (weak) sense, these vector fields are present
in the work of Schiffer \cite{Schiffer}}
\cite{Kir1}, \cite{KY}, also \cite{Kir2},
the Lie algebra $\Vect(S^1)$ acts 
on $\frK$ via  vector fields described in the following
way. Let $v(z)\frac\partial {\partial z}$ be a
real analytic  vector field
on the circle. Then the corresponding tangent vector 
at a point $f(z)\in\frK$ is 
\begin{equation}
\frac{f(z)^2}{2\pi i}
\int_{|t|=1}\frac{f'(t)^2 v(t)\,dt}
    {f(t)^2(f(t)-f(z))}
\label{eq:kirillov}
\end{equation}
The expressions of the vector fields $L_k$
in the terms of differential operators
in the coefficients $c_k$ is
\begin{align}
L_0&=\sum_{p>0}k c_k \frac\partial {\partial c_k}
\\
L_{k}&=\frac\partial{\partial c_k}+
\sum_{p>0}(1+p)c_k\frac\partial{\partial c_{k+p}},
\\
L_{-k}&=\sum_{p>0}
\Bigl\{\frac{1}{(2\pi i)^2}
\int\limits_{|h|=1-\epsilon}\int\limits_{|z|=\delta}
\frac{z^{-p-2} f^2(z) f'(h)^2 dh\,dz}
{h^{k-1}(f(h)-f(z))f^2(h)}\Bigr\}\cdot\frac\partial{\partial c_p}
\label{eq:l-minus-k}
=\\
&=\sum_{p>0}
\Bigl\{(k+p+1)c_{k+p}
+\frac{1}{(2\pi i)^2}
\int\limits_{|z|=1-\epsilon}\int\limits_{|h|=\delta}
\frac{z^{-p-2} f^2(z) f'(h)^2 dh\,dz}
{h^{k-1}(f(h)-f(z))f^2(h)}\Bigr\}\cdot\frac\partial{\partial c_p}
\label{eq:l-minus-k-1}
\end{align}
here $k>0$. 
The integral terms in (\ref{eq:l-minus-k}), (\ref{eq:l-minus-k-1}) 
differ by a position of a contour of integration;
the additional term in (\ref{eq:l-minus-k-1})  appears
since we move the contour through the pole $z=h$.

Also consider the series
\begin{equation}
\sum_{j=0}^\infty Q_n z^n
=z^2\left[h\frac{f'(z)^2}{f(z)^2}
+
\frac c{24} \Bigl(\frac{2 f'''(z)}{f'(z)}-
\frac {3 f''(z)^2}{f'(z)^2}\Bigl)\right]
\label{schwartz-generating}
\end{equation}
Then the operators (constructed in \cite{Ner-holomorphic}, 4.12)
\begin{equation}
\wh L_{k}:=L_{k},\qquad 
\wh L_{-k}:=L_{-k}+Q_k
\end{equation}
satisfy the relations
$$
[L_n,L_m]=(m-n)L_{m+n}+\frac c{12} (n^3-n)\delta_{n+m,0}
 $$
Also
$$
L_0 \cdot 1=h\cdot 1, \qquad L_{k}\cdot 1=0, \qquad
\text{for $k>0$}
$$
In other words, these operators determine a highest weight
representation of the Virasoro algebra
in the space of polynomials in variables $c_1$, $c_2$, \dots.
The corresponding action of the group of diffeomorphisms of the circle
is obtained in \cite{Ner-holomorphic}, 4.12.

\smallskip


{\bf \punct Some references.} The domain $\frK$
is an infinite-dimensional analogue of classical Cartan domains
and apparently the geometry and  harmonic analysis on Cartan domains
can be extended to $\frK$. Some elements
of this extension are known

\SS

-- K\"ahler structure, Ricci tensor, \cite{KY}, \cite{BR}

\SS

-- Reproducing kernels, \cite{Ner-holomorphic}, 4.9--4.10,
 in this paper they are
named as 'canonical cocycles' (also \cite{Ner-disser}) 

\SS

-- There is a theory of Verma modules, that are  dual to space
of holomorphic functionals on $\frK$, see
\cite{FF}, \cite{Ner-viniti}.

\SS

-- There are attempts to construct reasonable measures on
$\frK$ as limits of diffusions \cite{Mal}, \cite{AM}

\SS

-- Analogs of Olshanski semigroups for
$\frK$, \cite{Ner-semigroup}, \cite{Ner-holomorphic}.

\smallskip


{\bf\punct Purpose of the work.} 
It seems that the analysis
on $\frK$ requires some effective technique for 
manipulations with functionals in Taylor coefficients,
and we are trying to find some tools for this
(see also \cite{AR}, \cite{AB},\cite{Boi}).

We are trying to  write as explicitly as it is possible
the expressions for the action of the Virasoro algebra in terms
of Taylor coefficients, i.e.,
 (\ref{schwartz-generating}), (\ref{eq:l-minus-k}).
It seems for us that these expressions are simpler that it was
reasonable to wait (these expressions
are functions in infinite number of variables and they
can not be too simple), in any case factorial expressions
for coefficients with  single quadratic factors in
(\ref{eq:f-prime-f})--(\ref{eq:fra}), 
(\ref{eq:f-prime-prime-f})--(\ref{eq:frb}),
(\ref{eq:sch2})--(\ref{eq:sch3}),
(\ref{eq:L-k-1})--(\ref{eq:L-k-2})  were unexpected for us.

Also, we find a quasi-explicit expression for
the Grunsky  coefficients of univalent functions.

Our paper consists of manipulations with formal series, our
main tool is the Waring polynomials discussed in the Section 1;
but they are not present in final formulae in Section 3.
Calculations with Waring polynomials are contained in Section 2.


\section{Preliminaries. Waring and Faber polynomials}

\setcounter{equation}{0}

Here we recall some standard definitions and some 
well-known facts.

\smallskip


{\bf\punct\label{ss:waring-polynomials}
 Waring Polynomials.}
Let 
$$
a=\{a_j\}=(a_1,a_2,\dots)
$$
be an infinite sequence of complex numbers or of formal variables.
We define the {\it Waring polynomials} as $\cP_0(a)=1$, 
$$
\cP_n(a)=\cP_n(a_1,\dots,a_n), \qquad n=1,2,3,\dots
$$
by
\begin{equation}
\exp\Bigl\{
-\sum_{j=1}^\infty  \frac 1  j a_j z^j\Bigr\}
=\sum_{n=0}^\infty \cP_n(a) z^n
\label{eq:waring}
\end{equation}

Obviously,
$$
\cP_n(a_1,\dots,a_n)=\sum_{\mu_1,\mu_2,\dots,\mu_n:\,\sum j\mu_j=n}
 (-1)^{\sum\mu_j} \prod_j\frac{1}{j^{\mu_j} \mu_j!}
\prod_j a_j^{\mu_j}
$$


{\bf\punct \label{ss:faber} Faber polynomials.}
Let 
$$
h(z)=1+b_1 z+b_2 z^2+\dots
$$
We define the {\it Faber polynomials}
$$
\cQ_n(b)=\cQ_n(b_1,\dots, b_n), \qquad n=1,2,\dots 
$$
by
\begin{equation}
-\frac{zh'(z)}{h(z)}=-z\frac {d}{dz} \ln h(z)
=\sum_{n=1} \cQ_n(b) z^n
\label{eq:faber}
\end{equation}

Expanding $\ln h(z)=\ln(1+b_1z+b_2z^2+\dots)$ into a series, we obtain
$$
\frac 1n \cQ_n(b_1,\dots,b_n)=
\sum_{\mu_1,\dots,\mu_n:\,\sum j\mu_j=n}
(-1)^{\sum\mu_j} (\sum \mu_j-1)!
\prod \frac {{b_j}^{\mu_j}}{\mu_j!}
$$

\begin{proposition}
The maps
$$
b_1=\cP_1(a_1),\,\, b_2=\cP_2(a_1,a_2),\,\,b_3=\cP_3(b_1,b_2,b_3)
\,\,\dots
$$
and
$$
a_1=\cQ_1(b_1),\,\, a_2=\cQ_2(b_1,b_2),\,\,a_3=\cQ_3(b_1,b_2,b_3)
\,\,\dots
$$
are inverse.
\end{proposition}

{\sc Proof.} Let $\sum a_j z^j= -zh'(z)/h(z)$.
Then 
$$-\sum \frac{a_j}j z^j=\ln h(z)
\LR \exp\bigl\{-\sum \frac{a_j}j z^j\bigr\}=h(z)
$$      


{\bf\punct\label{ss:waring-th} Waring theorem.}
Let $x_\alpha$ be formal variables, $\alpha=1,2,\dots$. Denote
by 
$$p_n(x)=\sum_\alpha x_\alpha^n$$
 the Newton sums. By 
$$
e_n(x)=\sum_{i_1<i_2<\dots<i_n} x_{i_1}x_{i_2}\dots x_{i_n}
$$
we denote the elementary symmetric functions.

\begin{theorem}
$$
\P_n(-p_1(x),p_2(x)\dots,(-1)^{n}p_n(x))=e_n(x),
\qquad
\Q_n(e_1(x),e_2(x),\dots,e_n)=(-1)^{n}p_n(x)
$$
\label{th:waring}
\end{theorem}

{\sc Proof.} The first statement.
\begin{multline*}
\sum \cP_n(-p_1,\dots, (-1)^n p_n)z^n=
\exp\Bigl\{\sum\frac {(-1)^{j-1}}{j}
 (x_1^j+x_2^j+\dots)z^j\Bigr\}
=\\=
\prod_\alpha \exp\Bigl\{\sum_j\frac{(-1)^{j-1}}{j} (x_\alpha z)^j\Bigr\}
=\prod_\alpha (1+x_\alpha z)=
\sum_n e_n(x)z^n 
\end{multline*}

The second statement
\begin{multline*}
\sum_{k>0}\cQ_k(e_1,\dots,e_k)z^k=
-z\frac d{dz}\ln (1+e_1(x)z+e_2(x)z^2+\dots)
=\\=
-z\frac d{dz} \ln\prod_{j}(1+x_jz)
=-z\frac d{dz}\Bigl\{\sum \ln (1+x_jz)\Bigr\}
=\\=
-\sum_j \frac{zx_j}{1+x_jz}
=\sum_k (-1)^k z^k p_k(x)
\end{multline*}

{\sc Remark.} The Waring and the Faber polynomials and the Waring Theorem  appeared
in E.Waring's book "Meditationes algebraicae" 
1770. But it seems that the term "Waring polynomials"
is a neologism.
Waring polynomials appear in mathematics
in many situations, some ocassional
 list of references
is \cite{Che}, \cite{Sch}, \cite{Faa}
 \cite{Mac},  for instance,
characters of symmetric groups and inversion formula
for boson-fermion correspondence can be expressed in the terms
of these polynomials,\cite{Ner-boson}.

The term 'Faber polynomials' also has
a slightly different meaning, see the next subsection.


\bigskip

{\bf\punct Faber polynomials $\Phi_n(z)$.}
see \cite{Schiffer2}, see also \cite{Dur}.
Let 
\begin{align*}
h(z)&=1+b_1 z+b_2 z^2+\dots\\
g(z)&=zh(z^{-1})=z+b_1 +b_2z^{-1}+b_3 z^{-2}+\dots
\end{align*}
 The {\it Faber polynomial}
$\Phi_n(z)$ is a polynomial of one complex variable $z$
defined from the condition
\begin{equation}
\Phi_n\bigl( zh(z^{-1}))=
z^n+\sum_{k>0} \beta_{nk} z^{-k}
\label{eq:faber-usual}
\end{equation}

\begin{lemma}
$$
\Phi_n(z)=\Q_n(b_1-z,b_2,\dots,b_n)
$$
\end{lemma}

{\sc Proof,} see also \cite{Boi}. 
 Denote by $z=r(u)$ the function inverse to
$u=g(z)$ in a neighborhood of infinity. The equation
$$
z^n=
\Phi_n\bigl( zh(z^{-1}))-\sum_{k>0} \beta_{nk} z^{-k}
$$
is equivalent to 
$$
r(u)^n=\Phi_n(u)-\sum_{k>0} \beta_{nk} {r(u)}^{-k}
$$
Consider the Laurent expansion of $r(u)^n$,
$$
r(u)^n=\Bigl[u^n+p_{n-1}u^{n-1}+\dots + p_0 u^0\Bigr]
+ \Bigl\{p_{-1}u^{-1}+\dots\Bigr\}
$$
The expression in square brackets is $\Phi_n(u)$
and the expression in the curly brackets 
is $-\sum_{k>0} \beta_{nk} {r(u)}^{-k}$. Hence
$\Phi_n(u)$ is given by the Cauchy integral,
$$
\Phi_n(u)=\frac 1{2\pi i}\int_{|w|=R} \frac {r(w)^{n}dw}
 {w-u}
\qquad\text{where $|R|$ is large and $|u|<R$ also is large}
$$
After the substitution $w=g(z)$, we obtain
\begin{equation}
\Phi_n(u)=\int_C \frac{\zeta^n d(g(\zeta)-u)}{g(\zeta)-u}
\label{eq:faber-contour}
\end{equation}
Thus $\Phi_n(u)$ are the Laurent coefficients
of $\ln(\zeta h(\zeta^{-1})-u)$, and this easily implies
our statement.%
\QED
 
\SS


{\bf \punct Grunsky coefficients.} 
The {\it Grunsky coefficients} of a function
$$
g(z)=z+b_1 +b_2z^{-1}+b_3z^{-2}\dots
$$
are the numbers $\beta_{nk}$ from the formula
(\ref{eq:faber-usual}).

By (\ref{eq:faber-contour}),
$$
\sum_n\Phi_n(g(z))\zeta^{-n}=
\frac{\zeta g'(\zeta)}{g(\zeta)-g(z)}
$$
Transforming the left-hand side by (\ref{eq:faber-usual}),
we obtain
\begin{align}
\sum_{n=1}^\infty\sum_{k=1}^\infty
\beta_{nk} z^{-k}\zeta^{-n}=
\frac{\zeta g'(\zeta)}{g(\zeta)-g(z)}-
\frac{\zeta}{\zeta-z}
\\
=\zeta \frac\partial{\partial \zeta}
\ln\Bigl(\frac{g(\zeta)-g(z)}{\zeta-z}\Bigr)
\end{align}

The Grunsky matrix $\{\beta_{nk}\}$ has a fundamental
importance in the theory of univalent functions,
see \cite{Gru}, \cite{Gol}, \cite{Dur}, it also 
defines an $\Vect(S^1)$-equivariant embedding of
$\frK$ to infinite-dimensional Cartan domain
(observation of Yuriev) and appears as an element
in explicit formulae for representations
of the group of diffeomorphisms of circle
in Fock space \cite{Ner-holomorphic}, \cite{Ner-book}.


\section{Expansions of some products.}

\setcounter{equation}{0}


\smallskip

{\bf\punct} Let $\alpha_1$, $\alpha_2,\dots\in\C$,
$\alpha_j=0$ for big $j$. Let $\mu_j\in\C$.
Consider the weighted Newton sums
\begin{equation}
T_k(\alpha,\mu)=(-1)^k\sum \mu_j\alpha_j^k
\label{eq:w-newton}
\end{equation}

\begin{lemma}
\label{l:1}
\begin{equation}
\prod_j(1+\alpha_j s)^{\mu_j}=
\sum_{l\ge 0}\cP_l(T_1,T_2,\dots,T_l)s^l 
\label{eq:product-powers}
\end{equation}
\end{lemma}

{\sc Proof.}
\begin{multline*}
\prod_j(1+\alpha_j s)^{\mu_j}=
\exp\Bigl\{\sum_j\mu_j\ln(1+\alpha_j s)\Bigr\}
=\exp\Bigl\{
\sum_j \sum_k \frac{(-1)^{k-1}}k\mu_j\alpha_j^k s^k 
\Bigr\}
=\\=
\exp\Bigl\{-\sum \frac {1}{k} T_k(\alpha,\mu) s^k  \Bigr\}
\end{multline*}


{\bf\punct }
Let $\mu_2$,\dots, $\mu_3,\dots,\in\C$.
Let $\mu_k=0$ starting some number $k$.
Denote
\begin{equation}
N_p=N_p(\mu):=-\sum_{j\ge 2}
 \mu_j+\sum_{j\ge 2,\, j\,\, \text{divides}\,\, p} j\mu_j
\label{eq:Nl}
\end{equation}
 
\begin{lemma}
\label{l:2}
\begin{equation}
\prod_{j\ge 2}
\Bigl(\frac{1-t^j}{1-t}\Bigr)^{\mu_j}=
\sum_{l} \cP_l(N_1,N_2,\dots,N_l) t^l
\label{eq:N-l} 
\end{equation}
\end{lemma}

{\bf Proof}.
\begin{multline}
\prod_{j\geq 2}\Bigl(\frac{1-t^j}{1-t}\Bigr)^{\mu_j}=
\exp\Bigl\{\, -\,( \sum
\mu_j)\log(1-t)\, +\, \sum \mu_j\log(1-t^j)\, \Bigr\}
=\\=\exp\Bigl\{\, (\sum \mu_j)\sum_{k\geq 1}\frac{t^k}{k}\, -\, 
\sum_j \mu_j\sum_{m\geq 1}\frac{t^{jm}}{m}\Bigr\}
\label{A1}
\end{multline}
The coefficient in $t^k$ in the curly brackets
is 
$$
\frac 1k \sum \mu_j -  \sum_{\text{$j$ divides $k$}} \frac jk \mu_j
$$
and we transform (\ref{A1}) to
$$
\exp\Bigl\{\sum_{k\ge 1} \frac{-N_k(\mu)}k t^k\Bigr\}
$$
This implies the desired statement.%
\QED


\smallskip

{\bf\punct}
Let  $\bfm:=(m_1,m_2,\dots)$ be a sequence of
nonnegative integers, and $m_k=0$ for sufficiently
large $k$.
Denote
\begin{align*}
c^{\bfm}&:=c_1^{m_1} c_2^{m_2}\dots
\\
|\bfm|&:=m_1+m_2+\dots,\qquad
\|\bfm\|:=\sum j m_j
\end{align*}

\begin{lemma}
\label{l:3}
a) Let 
\begin{align*}
f(\xi)&=A_0+A_1\xi+A_2\xi^2+\dots\\
\theta(z)&=c_1z+c_2z^2+\dots
\end{align*}
 Then
\begin{equation}
f(\theta(z))=
\sum_\bfm A_{|\bfm|} \frac{|\bfm|!}{\prod m_j!} 
 c^\bfm z^{\|\bfm\|}
\label{eq:power-0}
\end{equation}

\begin{equation}
b)\quad\qquad (1+\theta(z))^p=
\sum_\bfm  \frac{p(p-1)\dots (p-|\bfm|+1)}{\prod m_j!} 
 c^\bfm z^{\|\bfm\|}
\qquad\qquad
\label{eq:power}
\end{equation}
\end{lemma}

{\sc Proof.}
 a) is easily deduced from multinomial
 identities, b) is a special case of a)
\hfill$\square$


\smallskip

{\bf\punct} Consider formal series
$\phi$, $\psi$ having the form
\begin{align*}
\phi(z)=1+c_1 z+ c_2 z^2+\dots\\ 
\psi(z)=\alpha_1 c_1z+\alpha_2 c_2 z^2+\dots
\end{align*}

\begin{lemma}
\label{l:4}
\label{l:psi-phi}
 Let $k$ be positive integer, $p\in \C$.
Then
\begin{equation}
\psi(z)^k \phi(z)^p=k!\sum_{\bfm}
\frac{p(p-1)\dots(p+k-|\bfm|+1)}
{\prod m_j!}
\cP_k(T_1,\dots, T_k) c^{\bfm} z^{\|\bfm\|}
\label{eq:psi-phi}
\end{equation}
where $T_k(\alpha)$ are given by
(\ref{eq:w-newton}).
\end{lemma}

{\sc Proof.}
Denote
\begin{equation}
U:=(s\psi(z)+\phi(z))^{k+p}=:
\sum H_j s^j
\label{eq:local}
\end{equation}
Assume additionally, that 
\begin{equation}
p\ne -1, -2,-3, \dots
\label{eq:artificial}
\end{equation}
 Then 
$$
 \psi(z)^k \phi(z)^p=\frac {k!}{(p+1)(p+2)\dots(p+k)} H_k 
$$
On the other hand,
$$
U=(1+(1+\alpha_1 s)c_1z+(1+\alpha_2s)c_2 z^2+\dots)^{k+p}
$$
We expand this expression in $z$ by (\ref{eq:power}),
$$ 
U=
\sum_{\bfm} \frac{(p+k)(p+k-1)\dots(p+k-|\bfm|+1)}
           {\prod m_j!}
\prod(1+\alpha_j s)^{m_j} \cdot c^{\bfm}z^{\|\bfm\|}
$$
Then we extract the coefficient in front of
$s^k$ by (\ref{eq:product-powers}).

Since final formula is continuous in $p$, we can omit the assumption
(\ref{eq:artificial}).
\QED


\smallskip

{\sc Remark.}  From (\ref{eq:product-powers}),
 if $k>|m|$, then ${\cal P}_k(T_1,  \cdots,T_k)=0 $.
If $k=|m|$, we understand  $p(p-1)\cdots(p+k-|m|+1)$ in (\ref{eq:psi-phi}) as $1$.

\smallskip

{\sc Remark.} Our calculations
allow to obtain the Taylor expansion of functions
$$
\psi(z)^k\phi(z)^p (\ln \phi(z))^m=
\frac {\partial^m} {\partial p^m} \psi(z)^k\phi(z)^p
$$
Actually we must differentiate only
the Pochhammer symbols in $p$.%
\QED


\smallskip

{\bf\punct} Now, let $\theta(z)=z+c_1z^2+\dots$,
Let $H(1+\xi)=A_0+A_1\xi+\dots$.

\begin{lemma}
\label{l:5}
$$
H\Bigl(\frac{\theta(z)-\theta(u)}{z-u}\Bigr)
=\sum_{\bfm}
A_{|\bfm|} \frac{|\bfm|! \,c^{\bfm}}{\prod m_j!}
\sum_{p,q: p+q=\|\bfm\|} 
\cP_p(\wt N_1,\dots,\wt N_p) z^p u^q
$$
where $\wt N_p$ are modified
$N_p$ from (\ref{eq:Nl})
$$
{\widetilde  N}_p=N_p(\mu)\quad with\quad \mu_j=m_{j-1}\quad for \, \, j\geq 2$$
\end{lemma}

{\sc Proof.} We must decompose
$$
A_0+\sum_{j>0} A_j\Bigl(c_1\frac{z^2-u^2}{z-u}   
+c_2\frac{z^3-u^3}{z-u}+\dots   \Bigr)^j
$$
The coefficient in $c^{|\bfm|}$ is 
$$
 A_{|\bfm|} \frac{|\bfm|!}{\prod m_j!}
\prod_{j\ge 1} \Bigl(\frac{z^{j+1}-u^{j+1}}{z-u}\Bigr)^{m_j}   
$$
and we apply Lemma \ref{l:2}.


\section{Formulae for univalent functions}

\setcounter{equation}{0}

{\bf \punct Schwarzian derivative.}
Let $\bfm=(m_1,m_2,\dots)$ be the same as above.
Denote
\begin{align*}
M_0&=|\bfm|=\sum m_j,\\
M_1&=\|\bfm\|=\sum j m_j\\
M_2&=\sum j^2 m_j
\end{align*}

\begin{proposition}
\label{pr:schwartz}
a) Let 
$$f(z)=z+c_1z^2+c_2 z^3+\dots$$
Then
\begin{equation}
\frac{z^{p+2} f'(z)^2}{f(z)^{p+2}}
=
\sum_{\bf m} \fra_{\bfm}(p) \frac{c^{\bfm}}{\prod m_j!}  z^{\|\bfm\|}
\label{eq:f-prime-f}
\end{equation}
where 
\begin{equation}
\fra_{\bfm}(p)=
(-1)^{|\bf m|}\frac{(p+|\bf m|-1)!}{(p+1)!}
\cdot
\Bigl\{p(p+1)+M_1^2-2(p+1)M_1-M_2\Bigr\}
\label{eq:fra}
\end{equation}
and $a_0(p)=1$

\begin{equation}
b)\qquad \frac{z^p f''(z)}{f(z)^p}=
\sum_{\bf m} \frb_{\bfm}(p) 
   \frac{c^{\bfm}}{\prod m_j!} z^{\|\bfm\|}
\qquad\qquad\qquad\qquad\qquad\qquad\qquad\qquad
\label{eq:f-prime-prime-f}
\end{equation}
where
\begin{equation}
\frb_{\bfm}(p)=\frac{(-1)^{|\bf m|+1}(p+|\bfm|-2)!}
{(p-1)! }\cdot\bigl\{M_1+M_2\bigr\} 
\label{eq:frb}
\end{equation}

c) For the Schwarzian derivative,
\begin{equation}
S_f(z)=2\frac{f'''(z)}{f'(z)}- 3 \frac{f''(z)^2}{f'(z)^2}
\label{eq:sch1}
\end{equation}
 we have the expansion
\begin{equation}
z^2 S_f(z)=\sum_{\bfm} \frd_{\bfm} 
\frac{ c^{\bfm}\prod (j+1)^{ m_j}}{\prod m_j!}\cdot z^{\|\bfm\|}
\label{eq:sch2}
\end{equation}
where
\begin{equation}
\frd_{\bfm}=(-1)^{|\bf m|} (|\bfm|-1)!\cdot\bigl\{M_2-3M_1^2+2M_1\}
\label{eq:sch3}
\end{equation}
\end{proposition}

In particular, the statement a) with $p=0$
and the statement c) gives explicit expressions for
the terms $Q_n$ in (\ref{schwartz-generating}).
 The coefficients when $p\neq 0$ in a) 
for any $p\in \Z$ and up to $||{\bf m}||=5$ have been explicited
in \cite{AM}, (A.1.7).

\smallskip

{\sc Proof.} a) We have the expression
$$
\frac{\bigl[1+(2c_1 z+3c_2z^2+\dots)\bigr]^2}{(1+c_1z+c_2z^2+\dots)^p}
$$
We apply Lemma \ref{l:psi-phi} 
to the functions $\phi(z):=f(z)/z=1+c_1z+c_2z^2+\dots$
and $\psi(z):=f'(z)-1=2c_1 z+3c_2z^2+\dots$

b) We have the expression
$$
z^{-1}\frac{1\cdot 2c_1z+2\cdot 3c_2z^2+ 3\cdot 4c_3z^3+\dots}
{(1+c_1z+c_2z^2+\dots)^p}
$$
 We apply Lemma \ref{l:psi-phi} to
$\phi(z)=f(z)/z$, $\psi(z)=zf''(z)$

c) follows from a),b). The factor $(1+j)^{m_j}$
appears since in our case $\phi(z)=f'(z)$.%
\QED

\smallskip


{\bf\punct The operators $L_{-k}$.}

\begin{proposition}
The vector field $L_{-k}$ at a point $f$ in $\frK$ is given by
\begin{equation}
\sum_{j=0}^\infty
\quad
\sum_{m:\,\|\bfm\|= j+k+1} 
\fra_{\bfm}(j+1) \frac{c^{\bfm}}
{\prod m_j!}
 f(z)^{j+2}  
\label{eq:L-k-1}
\end{equation}
or by
\begin{equation}
z^{1-k}f'(z)- \sum_{j=0}^{k} 
\quad
\sum_{m:\,\|\bfm\|= k-j} 
\fra_{\bfm}(-j)  
\frac{c^{\bfm}}
{\prod m_j!}
 f(z)^{1-j}  
\label{eq:L-k-2}
\end{equation}
where the coefficients $\fra$ are given (\ref{eq:fra})
\end{proposition}

{\sc Proof.}
First formula.
We must evaluate the integral
   (\ref{eq:kirillov}) with $v(t)=t^{-k+1}$, i.e.,
\begin{equation}
\frac 1{2\pi i}
\int_{|t|=1}
\frac{ f^2(z) f'(t)^2\,dt}
{t^{k-1}(f(t)-f(z))f^2(t)}
\label{eq:kirillov2}
\end{equation}
Emphasis that the integrand is
defined 
 for $0<|z|<1$,  $0<|t|<1$, $z\ne t$.
Since the expression is analytic in $z$,
we can think that $|z|$ is as small as we want.
We assume
$$|f(z)|<\min_{|t|=1}|f(t)|$$  
(since $f$ is univalent, by the Koebe--Bieberbach Theorem
the last minimum $\ge 1/4$,
see \cite{Gol}, \cite{Dur}).
Thus it is sufficient to
expand the expression
\begin{equation}
\frac{ f^2(z) f'(t)^2}
{(f(t)-f(z))f^2(t)}
\label{eq:expression-for-expansion}
\end{equation}
into a Laurent series 
$$
\sum U_k t^{k-2} 
$$
 in the domain
$|z|<|\epsilon|$, $1-\delta<|t|<1$ 
with small $\epsilon$, $\delta$.
The coefficients $U_k$ with positive $k$ are the vector fields
$L_{-k}$.
We transform (\ref{eq:expression-for-expansion}) to
$$
 f(z)^2f'(t)^2 f(t)^{-3}\Bigl(1+\frac{f(z)}{f(t)}+
\frac{f(z)^2}{f(t)^2}+\dots
 \Bigr)
=\sum_{j=0}^\infty \frac{f'(t)^2}{f(t)^{j+3}} f(z)^{j+2}
$$
and we apply Proposition \ref{pr:schwartz}.a.

\SS

The second formula. We fix $z\ne 0$ and
 move contour of integration in (\ref{eq:kirillov2})
to the position $|t|=\epsilon$. Passing through the pole
$t=z$, we obtain 
$$z^{1-k}f'(z)+
\frac 1{2\pi i}
\int_{|t|=\epsilon}
\frac{ f^2(z) f'(t)^2\,dt}
{t^{k-1}(f(t)-f(z))f^2(t)}
$$
To evaluate the integral, we look for the Laurent
expansion of the integrand
 in the domain $|t|< \epsilon$, $1-\delta<|z|<1$.
$$
- f(z)f'(t)^2f(t)^{-2}
\Bigl(1+\frac{f(t)}{f(z)}+
\frac{f(t)^2}{f(z)^2}+\dots
 \Bigr)
=
-\sum_{j=0}^\infty f'(t)^2 f(t)^{j-2} f(z)^{1-j}
$$
Now we apply Proposition \ref{pr:schwartz}.a.
Since we are interested only in 
 terms with
 powers of $t$ that are positive and less than $k$
the positive power of $t$, only finite number of terms
in (\ref{eq:L-k-2}) be present.%
\QED

\smallskip


{\bf\punct Grunsky coefficients.}

\begin{proposition}
$$
\beta_{nk}=n
 \sum_{\bfm:\, \|\bfm\|=n+k,\,m_1=0}
\frac {b^{\bfm} (|\bfm|-1)! }{\prod m_j!}
\cP_{k-|\bfm|}(L_1,\dots, L_{k-|\bfm|})
$$
where
$
L_k=N_k(m_3,m_4,\dots)
$
with $N$ given by (\ref{eq:Nl})
\end{proposition}

{\sc Proof.}  For $g(z)=z+b_1+b_2z^{-1}+ b_3 z^{-2}+\dots$,
we expand 
\begin{multline*}
\ln\frac{g(z)-g(\zeta)}{z-\zeta}=
\ln\Bigl(1+ 
b_2 \frac{z^{-1}-\zeta^{-1}}{z-\zeta}+
b_3 \frac{z^{-2}-\zeta^{-2}}{z-\zeta}+\dots
\Bigr)
=\\=
\ln\Bigl(1 -
b_2  z^{-1}\zeta^{-1} \frac{z^{-1}-\zeta^{-1}}{z^{-1}-\zeta^{-1}}-
b_3 z^{-1}\zeta^{-1}  \frac{z^{-2}-\zeta^{-2}}{z^{-1}-\zeta^{-1}}+\dots
\Bigr)
\end{multline*}

The coefficient of this expression in 
$b^{\bfm}=b_2^{m_2}b_3^{m_3}\dots$
is
\begin{multline*}
-\frac {(|\bfm|-1)! }{\prod m_j!} (z\zeta)^{-|\bfm|}
\prod_{j\ge 3}
\Bigl( \frac{z^{-(j-1)}-\zeta^{-(j-1)}}
{z^{-1}-\zeta^{-1}}\Bigr)^{m_j}
=\\=
-\frac {(|\bfm|-1)! }{ \prod m_j!} (z\zeta)^{-|\bfm|}
\sum_{j\ge 0} \cP_j(L_1,\dots,L_j) z^{-j} \zeta^{-(\|\bfm\|-2|\bfm|-j)}
\end{multline*}
we take $n=||\bf m||-|\bf m|-j$, $\, k=|\bf m|+j$ and this is our statement.%
\QED

\bigskip


{\bf Acknowledgements.} This work is issued from discussions
while Yu. Neretin visited the University of Picardie Jules Verne
and  the LAMFA CNRS UMR 6140 (Amiens) in February 2006.

 \tt H.Airault:

INSSET, Universite de Picardie Jules Verne, 48, rue Raspail, 02100
Saint-Quentin(Aisne)
Laboratoire CNRS UMR 6140, LAMFA, 33 rue Saint-Leu, 80039 Amiens, France.

email: hairault@insset.u-picardie.fr

\smallskip

 Yu. Neretin:

Department of Mathematics

 University of Vienna, Nordbergstrasse, 15, Vienna, Austria

\&

 Institute for Theoretical and Experimental Physics,

B. Cheremushkinskaya 25, Moscow, Russia 117259.

email: neretin@mccme.ru

URL:www.mat.univie.ac.at/$\sim$neretin

\end{document}